\documentclass{article}
\usepackage{amssymb}
\usepackage{indentfirst}

\makeatletter

\@addtoreset{equation}{section}
\makeatother

\def\vs{\vskip.5cm}

\def\P{{I\kern-0.3emP}}
\def\Q{{\P_A}}
\def\pa{\P(A)}
\def\fa{f_A}

\def\za{\zeta_A}

\def\ta{\tau_A}    
\def\ca{ r_A }     
\def\na{ n_A }     
\def\oa{ \tau(A) } 

\def\del{ \delta }

\def\eps{ \epsilon }
\def\sd{\overline{d}}

\def\lima{\rho(A)}

\def\Nset{I\kern-0.3emN}
\def\Zset{{Z\kern-0.6emZ}}
\def\Rset{I\kern-0.3emR}
\def\bbe{I\kern-0.3emE}
\def\one{1\kern-0.3em1}

\def\calc{\mathcal C}

\def\calr{\mathcal R}
\def\ra{ \calr(A) }
\def\ps{ \mathrm }

\def\lp{ \left( }
\def\rp{ \right) }
\def\ll{ \left\{ }
\def\rl{ \right\} }
\def\lc{ \left[ }
\def\rc{ \right] }
\def\lv{ \left\vert }
\def\rv{ \right\vert }

\def\go{ \rightarrow }
\def\Go{ \Rightarrow }
\def\Goback{ \Leftrightarrow }

\def\pros{ $\ll X_m \rl_{m \in \Zset}$ }
\def\amix{ $\alpha$-mixing }
\def\fmix{ $\phi$-mixing }
\def\pmix{ $\psi$-mixing }

\def\beq{ \begin{displaymath} }
\def\eeq{ \end{displaymath}  }
\def\beqn{ \begin{equation} }
\def\eeqn{ \end{equation}  }
\def\beqa{ \begin{eqnarray*} }
\def\eeqa{ \end{eqnarray*} }
\def\beqan{ \begin{eqnarray} }
\def\eeqan{ \end{eqnarray} }
\def\nn{ \nonumber }

\newtheorem{teo}{Theorem}[section]
\newtheorem{defi}{Definition}[section]
\newtheorem{eje}{Example}[section]
\newtheorem{lem}{Lemma}[section]
\newtheorem{pro}{Proposition}[section]
\newtheorem{cor}{Corollary}[section]
\newtheorem{rem}{Remark}[section]

\newtheorem{pf*}{Proof}

\def\bteo{ \begin{teo} }
\def\eteo{ \end{teo} }
\def\bdefi{\begin{defi} }
\def\edefi{ \end{defi} }
\def\beje{ \begin{eje} }
\def\eeje{ \end{eje} }
\def\blem{ \begin{lem} }
\def\elem{ \end{lem} }
\def\bpro{ \begin{pro} }
\def\epro{ \end{pro} }
\def\bcor{ \begin{cor} }
\def\ecor{ \end{cor} }
\def\brem{ \begin{rem} }
\def\erem{ \end{rem} }

\def\qed{\Box}
\def\bpf{ \noindent{\bf Proof } }
\def\epf{ $\qed$ \vs}

\title{ Sharp error terms for return time statistics under mixing conditions
}
\author{Miguel Abadi
\thanks{ IMECC, Universidade Estadual de Campinas,
         P\c ca S\'ergio Buarque de Holanda 651 Cid. Univ.
         CP 6065, Cep. 13083-859, Campinas SP, Brazil. Phone +55-19 37886021 
         \rm{miguel@ime.unicamp.br}  }
\and Nicolas Vergne
\thanks{ Universit\'{e} d'Evry Val d'Essonne, D\'{e}partementcMath\'{e}matiques,
         Laboratoire Statistique et G\'{e}nome, 91 000 Evry, France.
         \rm{vergne@genopole.cnrs.fr}  }
}
\date{}

\begin{document}

\maketitle

\begin{abstract}
We describe the statistics of repetition times of a string of
symbols in a stochastic process.

Denote by $\ta$ the time elapsed until the process spells the finite string $A$
and by $S_A$ the number of consecutive repetitions of $A$.
We prove that, if the length of the string grows unbondedly, 
(1) the distribution of $\ta$, when the process starts with $A$, is well approximated
    by a certain mixture of the point measure at the origin and an exponential law, and 
(2) $S_A$ is approximately geometrically distributed. 
We provide sharp error terms for each of these approximations.
The errors we obtain are point-wise and allow to get also approximations for all the
moments of $\ta$ and $S_A$.
To obtain (1) we assume that the process is $\phi$-mixing while to obtain (2) we assume 
the convergence of certain contidional probabilities.
\end{abstract}
{\bf Running head:} Return times under mixing conditions.

\section{\bf INTRODUCTION}


This paper describes the return time statistics  of a string
of symbols in a mixing stochastic process with a finite alphabet.
Generally speaking, the study of the time elapsed until the
first occurrence of a small probability event has a long history. See for instance \cite{GS} and the
references therein.
The typical result is:
\beqn \label{eqtyp}
\lim_{n \go \infty} \P\lp \tau_{A_n} > t \ b_n \ | \ \mu_0 \rp = e^{-t} \ .
\eeqn
Here  $\tau_{A_n}$ is the first time the process hits a given
measurable set $A_n$, $n\in\Nset$. We assume that the measure
$\P \lp A_n \rp$ goes to zero as $n \go \infty$.
$\{b_n\}_{n \in \Nset}$ is a suitable
re-scaling sequence of positive numbers and $\mu_0$
is a given initial condition.

Recently an exhaustive analysis of these statistics
was motivated by applications in  different
areas as entropy estimation, genome analysis, computer science, linguistic, among others.
From the point of view of applications, a fundamental task is to understand 
the rate of convergence of the limit (\ref{eqtyp}).
A detailed review of such results appearing in the literature 
can be found  in \cite{AG}. 

It is the purpose of this paper to present the following new results: 
For \emph{any} string $A$ of lenght $n$\\
- A sharp upper bound for the above rate of convergence that holds when $\mu_0=A$. 
  In this case we assume that the process is \fmix.\\
- A sharp upper bound for the difference between the law of the number of 
  \emph{consecutive} visits to $A$ and a geometric law. In this case we assume a kind of continuity for 
  certain conditional probabilities, see (\ref{eq:lim}).
  on the 

\vs
When $\mu_0$ is taken as $A$, we refer to the distribution $\P(\ta > t \ | \ A)$
as the \emph{return time}.
In general it can not
be well approximated by an exponential law.
This was firstly noted by Hirata,
when he proved the convergence of the number of visits to a small cylinder
around a point to the Poisson law.
His result holds for axiom A diffeomorphisms (see \cite{hir}).
The result holds for \emph{almost} every point.
Then, he proved that for periodic points,
the asymptotic limit law of the return time
differs from the one-level Poisson law, namely $e^{- t}$.

Our first result concerns the rate of convergence of
limit in (\ref{eqtyp}) when $\mu_0=A$ for any string $A$ of lenght $n$..
We prove that the return time law
converges to a 
convex combination of a Dirac law at the origin and an exponential law.
Specifically, we show that for large $n$
\beq 
\P\lp \ta > {t\over\pa} \ | \  A \rp  \approx \ll
\begin{array}{ll}
1 & t\le \pa\oa \\
\za e^{-\za t}& t>\pa\oa  
\end{array}
\right. \ .
\eeq 
$\oa$ is the position of the first overlap of $A$ with a copy of itself 
(see definition below).
$\za$ is a parameter related to the overlap properties of the string $A$.
It is worth noting that the parameter of the exponential law is exactly
the weight of the convex combination.
So far, the overlap properties of a string appears as a major factor
to describe the statistical properties of the return time.
For instance, if a string overlaps itself, then
it will turn out in the sequel that $\za\not=1$ and the return time distribution
approximates the above mixture of laws.
However, for a word which does not overlap itself,
it will turn out that $\za=1$ and the return time distribution approximates
a purely exponential law.
For the role of overlaps an a treatment of the independent case
with a good introduction to the previous literature see \cite{BT}, and for
the Markov case with a probability generating functions point of view see \cite{ste}.

It is worth recalling at this point that when in equation (\ref{eqtyp}) the initial condition
is the equilibrium measure of the process,
$\ta$ is called the {\em hitting time} of $A$.
In \cite{HSV} it is proved
a rate of convergence of the return time
as function of the distance between the hitting time
and return time laws. While this result 
applies only for cylinders around non-periodic points,
our result applies  to {\em all} of them.

The great enhancement of our work is that,
contrarily to all the previous works which present bounds depending only on the string $A$,
our error estimate decays exponentially fast in $t$
for all  $t>0$.
As a byproduct we obtain \emph{explicit expressions}  for \emph{ all the moments}
of the return time.
This also appears as a generalization of
the famous Kac's lemma (see \cite{kac}) which states
that the \emph{first} moment of the return time to a string $A$ of positive measure
is equal to $\pa^{-1}$ and the result in
\cite{chaz} which presents
conditions for the \emph{existence} of the moments of
return times.
Further, \cite{HSV} proves that hitting and return times
coincide if and only if the return time converges to the exponential law.
We extend this result establishing that
the laws of hitting and return times coincide if and only if
the weight of the Dirac measure in the convex combination of the return time
law is zero, which is equivalent to consider a non-overlapping string.

Our framework is the class of \fmix processes.
For instance, irreducible and aperiodic finite state
Markov chains are known to be $\psi$-mixing (and then \fmix) with exponential decay.
Moreover, Gibbs states which have  summable variations are $\psi$-mixing (see \cite{wal}).
They  have exponential decay if they have
H\"older continuous potential (see \cite{bow}).
However, sometimes the \pmix condition is very restricted hypothesis
difficult to test. 
We establish our result under the more general \fmix condition.
Further examples of \fmix processes can be found in \cite{LSV}.
The error term is explicitly expressed as a function of the mixing rate $\phi$.
We refer the reader to \cite{dou}  for a source of examples and definitions 
of the several kinds of mixing processes.

The base of our proof is a sharp upper bound on the
rate of convergence of the hitting time
to an exponential law proved in \cite{aba4}.

\vs
The self-repeating phenomena in the distribution
of the return time leads us to consider the problem of
the sojourn time.
Our second  result states that the law of the number of consecutive
repetitions of the string $A$, denoted by $S_A$,
converges to a geometric law. Namely
\beqn \label{typg}
\P\lp S_A = k \ | \ A \rp \approx (1-\rho(A))\rho(A)^{k} \ .
\eeqn
Again here, the parameter $\rho(A)$ depends on the overlap properties
of the string.
Furthermore we show that under suitable
conditions one has $\rho(A) \approx 1-\za$.
As far as we know, this is the first result on
this subject for dependent processes.

As in our previous result, 
the error bound we obtain decreases geometrically fast in
$k$ (see (\ref{typg})).
This decay on the error bound  allows us to obtain
an approximation for {\em all the moments} of $S_A$ for those of
a geometrically distributed random variable.
\vs
Our results are applied in a forthcoming paper:
In \cite{AHV} the authors prove large deviations
and fluctuations properties of the repetition time function
introduced by Wyner and Ziv in \cite{WZ} and further by
Ornstein and Weiss in \cite{OW}, and get entropy estimators.
\vs
This paper is organized as follows.
In section 2 we establish our framework.
In section 3 we describe the self-repeating
properties needed to state the return time result.
In section 4 we establish the approximation
for the return time law. This is Theorem \ref{teo:ret}.
Finally, in section 5 we state and prove the geometric approximation for the
consecutive repetitions  of a string. This is Theorem \ref{teo:geo}.

\section{\bf FRAMEWORK AND NOTATION}

Let ${\cal C}$ be a finite set. Put $\Omega={\mathcal C}^Z$.
For each $x=(x_m)_{m\in Z} \in\Omega$ and
$m \in Z,$ let $X_m:\Omega \rightarrow {\mathcal C}$ be
the $m$-th coordinate projection, that is $X_m(x)=x_m$.
We denote by $T:\Omega \rightarrow \Omega$
the one-step-left shift operator, namely $(T(x))_m= x_{m+1}$.

We denote by ${\mathcal F}$ the $\sigma$-algebra over $\Omega$
generated by strings. Moreover
we denote by ${\mathcal F}_{I}$ the $ \sigma$-algebra generated
by strings with coordinates in  $I$, $I \subseteq Z$.

For a subset $A \subseteq \Omega $, $A \in {\mathcal C}^n$
if and only if
\[
A = \{ X_0 = a_0; \dots ; X_{n-1} = a_{n-1} \} \ ,
\]
with $a_i \in {\mathcal C}, \; i=0, \dots ,n-1$.

We consider an invariant probability measure $\P$ over ${\mathcal F}$.
We shall assume without loss of generality that there is no singleton  of
probability 0.

For two measurable sets  $V$ and $W$, we denote as usual
$\P\lp V  |  W \rp= \P_W(V) = \P\lp V  ;  W \rp / \P(W)$
the conditional measure of $V$ given $W$.
We write
$\P\lp V ; W \rp=\P\lp V \cap   W\rp$.

We say that the process \pros is
\fmix if the sequence
\[
\phi(l)=
\sup \lv \P_B(C)  - \P(C) \rv \ ,
\]
converges to zero.
The supremum is taken over $B$ and $C$ such that
$B \in {\mathcal F}_{\{0,.,n\}}, n \in \Nset,\P(B)>0, C \in {\mathcal F}_{\{m\in\Nset | \ m\ge n+l+1 \}}$.

We use the measure theoretic notation:
$\{ X_{n}^{m}=x_{n}^{m} \} = \{X_n=x_n, \dots , X_m=x_m \}. $
For an $n$-string $A=\{ X_{0}^{n-1}=x_{0}^{n-1} \} $ and $1\le w \le n$,
we write $A^{(w)}= \{ X_{n-w}^{n-1}=x_{n-w}^{n-1} \} $ for the $w$-string
belonging to the $\sigma$-algebra ${\mathcal F}_{\{n-w,\dots,n-1\}}$
and consisting of the {\em last} $w$ symbols of $A$.
We write $V^c=\Omega\backslash V$, for
the complement of $V$.

The conditional mean of a r.v. $X$ with respect to any measurable set $V$ will be
denoted by $\bbe_V(X)$ and we put $\bbe(X)$ when $V=\Omega$.
Wherever
it is not ambiguous we will write $C$ for different positive
constants even in the same sequence of equalities/inequalities.
For brevity we put $(a \lor b)=\max\{a,b\}$ and
$(a \land b)=\min\{a,b\}$.

\section{PERIODS}

\bdefi
Let $A \in {\mathcal C}_n$. We define  the {\em period} of  $A$
(with respect to $T$) as the number $\oa$ defined as follows:
\[
\oa=\min\ll k \in \{1,\dots,n\} \; | \; A\cap\ T^{-k}(A)\ne\emptyset \rl \ .
\]
\edefi

By definition, if $A\in\calc_n$, then $A=(a_0,\dots,a_{n-1}), a_i\in \calc$ for $0\le i \le n-1$.
For instance, pick up  $A= (aaaabbaaaabbaaa) \in \calc^{15}$. Then shift a copy of $A$ until there is a fit between them.
Namely
\beq 
\begin{array}{rll}
       A= & {\tt aaaabb}&\!\!\!\!\! {\tt aaaabbaaa }\\
T^{-6}(A)=&      & \!\!\!\!\! {\tt aaaabb  aaa  \ a  bbaaa}\\
\end{array}
\ .
\eeq
Notice that there is no fit between $A$ and $T^{-j}(A)$
if $j=1,\dots,5$. So that $\oa=6$.

Let us take $A \in {\mathcal C}_n$, and write $n=q \, \oa + r$, with $q=[n/\oa]$ and
$0\le r < \oa$. Thus
\[
A
= \ll X_{0}^{\oa-1}=X_{\oa}^{2\oa-1}= \dots
= X_{(q-1)\oa}^{q\oa-1}=a_{0}^{\oa-1} \ ; \ X_{q\oa}^{n-1}=a_{0}^{r-1} \rl \ .
\]
So, we say that $A$ has  {\em period} $\oa$ and  {\em rest }$r$.
We remark that periods can be ``read backward''
(and for the purpose of section 5 it will be more useful to do it in this way),
that is
\beqa
A \!\!\!\!\!
&=& \!\!\!\!
\ll X_{0}^{r-1}=a_{n-r}^{n-1}  ; 
X_{n-q\oa}^{n-(q-1)\oa-1}=..=X_{n-2\oa}^{n-\oa-1}=X_{n-\oa}^{n-1}=a_{n-\oa}^{n-1}\rl \\
&=& \!\!\!
 \bigcap_{j=1}^{(q-1)\oa} T^{j\oa}(A^{(\oa)}) \ \cap \  T^{q\oa}(A^{(r)})  \ .
\eeqa
We recall the definition of $A^{(w)}, 1\le w\le n$, from the end of section 2.
For instance, using the previously chosen $A$,
\beqn \label{string}
A= (
\overbrace{\tt aaaabb}^{ \ps{ period} }
\overbrace{\tt aaaabb}^{ \ps{ period} }
\overbrace{\tt aaa}^{ \ps{ rest} }
)= (
\overbrace{ \underbrace{\tt aaa}_{    T^{12}A^{(3)} }  }^{ \ps{ rest} }
\overbrace{ \underbrace{\tt abbaaa}_{ T^{6}A^{(6)}  }  }^{ \ps{ period} }
\overbrace{ \underbrace{\tt abbaaa}_{ A^{(6)}  }  }^{ \ps{ period} }
) \ .
\eeqn
In the middle of the above equality, periods
are read forward while in the right hand side
periods are read backward.

\vs

Consider the set of overlapping positions of $A$:
\[
\ll k \in \{1,\dots,n-1\} \ | \ A \cap \  T^{-k}(A) \ne \emptyset
\rl = \{ \oa,\dots, [n/\oa]\oa \} \cup \ra \ ,
\]
where
\[
\ra = \ll k \in \{ [n/\oa]\oa+1,\dots,n-1 \}  \ | \ A \cap \ T^{-k}(A)
\ne \emptyset \rl\ .
\]
The set  $ \{ \oa,\dots, [n/\oa]\oa \}$ is called the set of principal periods
of $A$ while $\ra$ is called the set of secondary periods of $A$.
Furthermore, put $\ca = \# \ra$.
Observe that one has  $ 0 \le \ca  <n/2$.

\vs
The notion of period is related to the notion of \emph{retun times}.
\bdefi
Given $A \in {\mathcal C}_n$, we define the \emph{hitting time}
${\ta}: \Omega \rightarrow \Nset \cup \{\infty\}$ as the following
random variable:
For any $x \in \Omega$
\[
\ta(x) = \inf\{k \ge 1 : T^k(x) \in A \} \ .
\]
The \emph{return time} is the hitting time restricted to the set $A$,
namely $\ta|_A$. 
\edefi

We remark the difference between $\ta$ and $\oa$: 
while $\ta(x)$ is the first time $A$ appears in $x$, $\oa$ is the first overlapping
position of $A$.

Return times before $\oa$ are not possible, thus, $\Q\lp \ta < \oa \rp=0$.
Still, if $A$
does not return at time $\oa$, then it can not return at times $k
\oa$, with $2\le k \le [n/\oa]$, so one has
\[
\Q\lp  \oa < \ta  \le [n/\oa] \oa \rp=0 .
\]
The first possible return time after $\oa$ is
\beq
\na = \ll  
\begin{array}{lr}
\min\ra  & \qquad    \ra\ne\emptyset \\
n    & \qquad    \ra=\emptyset 
\end{array} \right. \ .
\eeq
Furthermore, by definition of $\ra$ one has
$ A \bigcap T^{-j}(A)=\emptyset $ 
for all $j$ such that $[n/\oa] \oa <  j  \le n-1 $ and $j\not\in\ra$.
Thus
\[
\Q\lp \ll [n/\oa] \oa +1 \le  \ta  \le n-1 \rl \cap \ll  \ta\not\in\ra \rl \rp=0 .
\]
We finally  remark that
$$ T^{-i}A \cap T^{-j}A =\emptyset \qquad \forall i,j \in\ra \ .$$ 
Otherwise
it would contradict the fact that the first return time to $A$ is $\oa$ since
for $i,j\in\ra$ one has $|i-j|< \oa $.
We conclude that
\beqn \label{:ra}
\Q\lp T^{-i}A \cap T^{-j}A \ | \ i,j\in \ra  \rp=0 .
\eeqn

\section{RETURN TIMES}

For $A \in {\mathcal C}_n$ define
\[
\za \stackrel{def}{=}
 \Q(\ta\ne \oa)=\Q(\ta > \oa ) \ .
\]
The equality follows by the comment at the end of the previous section.

It would be useful for the reader to note now that according to the comments
of the previous section, one has
\beqn \label{imagen}
\ta|_A \ \in \ \{\oa\} \cup \ra \cup \{ k \in \Nset \ | \ k \ge n \} \ .
\eeqn
 
We now introduce the error terms that appear in the statement of our main result 
of this section.

\bdefi Let us define
%
\beqn \label{:ret2}
\eps(A)
\stackrel{def}{=}
\inf_{0\le w \le \na} \lc (2n+\oa) \P(A^{(w)}) + \phi\lp \na - w \rp\rc  \ .
\eeqn
\edefi

\bteo \label{teo:ret}
Let \pros be  a \fmix process.
Then, 
for all $A \in \calc^n, n\in \Nset$
the following inequality holds for all $t$:
\beqn
\lv \Q \lp \ta > t \rp
- \one_{\{t< \oa\}} - \one_{\{t\ge\oa\}}\za e^{- \za\pa (t-\oa)} \rv
\le 54  \eps(A)  f(A,t)  , \label{eqteo:ret}
\eeqn
where $f(A,t)= \pa t  e^{- ( \za-16\eps(A) ) \pa t}$.
\eteo

We postpone an example showing the sharpness of $\eps(A)$ after
Lemma \ref{lem:b1}.

\brem
$A^{(n_A)}$ is the part of the string $A$ which does not overlap
itself in $A\cap T^{-n_A}A$.
Note that $n_A$ is the position of the first possible return time 
after $\oa$.
Recall that $r_A=\# \ra$ and $n_A=n$ if $\ra=\emptyset$.
Thus $A^{(w)}$ with $1\le w \le \na$ is the part of the string $A^{(n_A)}$ 
after taking out its first $n_A-w$ letters 
(this will be to create a gap of length $n_A-w$ to use the mixing property).
\erem

\brem
When $\ra=\emptyset$, namely, $A$ does not have secondary periods, 
the error $\eps(A)$ of Theorem \ref{teo:ret} becomes 
$\inf_{0\le w \le n} \lc n \P(A^{(w)}) + \phi(n-w)\rc \ . $
\erem

\brem
In the error term of the theorem, $\eps(A)$ provides a bound which shows the
convergence uniform in $t$
of the return time law to that mixture of laws as
the length of the string growths.
The factor $\pa t$ provides an extra bound for values of $t$ smaller than $1/\pa$.
The factor $e^{-(\za-16\eps(A)))\pa t}$ provides an extra bound for values of $t$ larger than $1/\pa$.
\erem

\brem On one hand $\pa\le Ce^{-c n}$ (see \cite{aba1}). 
On the other hand, by construction $\na > n/2$. Further
$\phi(n)\go 0$ as $n\go\infty$. Taking for instance $w=n/4$ in (\ref{:ret2})
we warrant the smallness of $\eps(A)$ for large enough $n$.
\erem

\bcor \label{:kac} Let the process \pros be $\phi$-mixing. 
Let $\beta>0$. Then, for all $A \in \calc^n, n\in \Nset$, the
$\beta$-moment of the re-scaled time $\pa\ta$ approaches, as
$n\go\infty$, to $\Gamma(\beta+1)/\za^{\beta-1}$. Moreover
\beqn \label{eq:kac}
\lv \pa^\beta \bbe_A ( \ta^\beta )  -  \frac{\Gamma(\beta+1)}{ \za^{\beta-1} }\rv
\le  \eps^*(A)  {C \beta \ e^{2\eps(A)(\beta+1)/\za} \over\za^2} \  {\Gamma(\beta+1)\over \za^{\beta-1} }\ ,
\eeqn
where $\eps^*(A)= (\eps(A)\lor (n\pa)^\beta) $, $C>0$ is a constant 
and $\Gamma$ is the analytic gamma function.
\ecor

\brem
In particular, the corollary establishes that all the moments of the return time are finite.
\erem

\brem
In the special case when $\beta=1$, the
above corollary establishes a weak version of Kac's Lemma (see \cite{kac}).
\erem

\brem
For each $\beta$ fixed and $n$ large enough one has 
${\beta \ e^{2\eps(A)(\beta+1)}/\za^2}$ is close to $\beta/\za^2$.
Thus in virtue of inequality (\ref{eq:kac}), the corollary  reads not just as a difference 
result but also as a ratio result.
\erem

The next corollary extends Theorem 2.1 in
\cite{HSV}.

\bcor \label{cor:equiv}
Let the process \pros be $\phi$-mixing.
There exists a constant $C>0$ such that, for all
$A\in \mathcal{C}^n, n\in \Nset$ and all $t>0$ 
the following conditions are equivalent:
\begin{enumerate}
\item[(a)]
$
\lv \Q \lp \ta > t \rp - e^{-\pa t} \rv  \le  C \ \eps(A) \  f(A,t)\ ,
$
\item[(b)]
$
\lv \Q \lp \ta > t \rp -\P \lp\ta > t \rp \rv \le C \ \eps(A) \ f(A,t) \ ,
$
\item[(c)]
$
\lv \P \lp \ta > t \rp - e^{-\pa t} \rv  \le  C \ \eps(A) \  f(A,t)\ ,
$
\item[(d)]
$
\lv \za  - 1 \rv \leq C \ \eps(A)   \ .
$
\end{enumerate}
Moreover, if $\{A_n\}_{n\in\Nset}$ is a sequence of strings
such that $\P(A_n)\go 0$ as $n\go \infty$,
then the following
conditions are equivalent:

\noindent ($\tilde{a}$) the return time law of $A_n$
converges to a parameter one exponential law,\\
\noindent ($\tilde{b}$) the return time law
and the hitting time law of $A_n$ converge to the same law,\\
\noindent ($\tilde{c}$) the hitting time law of $A_n$
converges to a parameter one exponential law,\\
\noindent ($\tilde{d}$) The sequence $(\zeta_{A_n})_{n\in\Nset}$ converges to one.
\ecor

\subsection{Preparatory results}

Here we collect a number of results that will be useful for the proof of Theorem \ref{teo:ret}.
In what follows and for shorthand notation we put $\fa=1/(2\pa)$ (factor 2 is rather technical).
The next lemma is a useful way to use the \fmix property.

\blem \label{fmix} Let \pros be  a \fmix process.
Suppose that $A \supseteq B \in {\mathcal F}_{\{0,\dots,b\}},
C \in {\mathcal F}_{\{x\in\Nset | \ x\ge b+n\}}$ with $b,g \in \Nset$.
The following inequality holds:
\[
\Q( B ; C ) \le \Q(B) \lp \P(C) + \phi(n) \rp \ .
\]
\elem
\bpf
Since $B \subseteq A $, obviously
$\P( A \cap B \cap C ) = \P( B \cap C )$. By the \fmix
property
$
 \P( B ; C ) \le \P(B) \lp \P(C) + \phi(n) \rp  .
$
Dividing the above inequality by $\pa$ the lemma follows.
\epf

The following lemma says that return times over $\ra$ have small probability.

\blem \label{lem:b1}
Let \pros be a \fmix process.
For all $A \in \calc^n$,
the following inequality holds:
\beqn \label{eq:b1}
\Q\lp \ta\in\ra \rp
\le  \eps(A) \ . 
\eeqn
\elem

\bpf
For any $w$ such that $1\le w \le n_A$
\beqan
     \Q\lp \ta\in\ra \rp
&\le&  \Q\lp \bigcup_{j\in \ra} T^{-j} A       \rp \nn\\
&\le&\Q\lp \bigcup_{j\in \ra} T^{-j} A^{(w)} \rp \nn\\
&\le& \ca\P\lp A^{(w)} \rp+ \phi( \na-w )  \ . \label{model}
\eeqan
The first inequality follows by (\ref{:ra}).
Since $T^{-j}A \subset T^{-j}A^{(w)}$,
second one follows.
Third one follows by the above lemma
with $B=A$ and $C=\cup_{j\in \ra} T^{-j} A^{(w)} $.
This ends the proof since $w$ is arbitrary.
\epf

\beje
Consider a process \pros defined on the alphabet
${\mathcal C}=\{a,b\}$. Consider the string introduced in (\ref{string}):
\[A=\ll (X_0...X_{14})=( \tt aaaabbaaaabbaaa) \rl \ .\]
Then, $n=15$, $\oa=6$, $\ra=\{13,14\}$, $\ca = 2$ and $\na = 13$.
Thus 
\[A^{(13)}=\ll  (X_2...X_{14})=( \tt aabbaaaabbaaa) \rl \ .\]
The \fmix property factorizes the probability
\beqa
\Q\lp \bigcup_{j=13}^{14} T^{-j} A \rp 
=
\Q\lp \bigcup_{j=13}^{14} T^{-j} A^{(13)} \rp 
\le
\Q\lp \bigcup_{j=13}^{14} T^{-j} A^{(w)} \rp \ .
\eeqa
In such case, a gap at $t=15$ of length $w$ with
$0\le w \le 13 $ is the best we can do to apply the \fmix property.
\eeje

The next lemma will be used to get the non-uniform factor $f(A,t)$ in the
error term of Theorem \ref{teo:ret}.

\blem \label{itera}
Let  \pros be a  \fmix  process.
Let $A \in \calc^n$ and
let $B \in \mathcal{F}_{\{x\in\Nset | \ x\ge k\fa\}}$, with $k \in \Nset$.
Then the following inequality holds:
\beq 
\Q\lp \ta > k \fa \ ; \ B \rp 
\le
\lc \P\lp \ta > \fa-2n\rp + \phi(n) \rc^{k-1}
\lc \P(B) + \phi(n)  \rc \ .
\eeq 
\elem
\bpf
First introduce a gap of length $2n$ between $\{\ta > k \fa \}$ and $B$.
Then use Lemma \ref{fmix} to get the inequalities
\beqan \label{1itera}
\Q( \ta > k \fa \ ; \ B )
&\le&
\Q(\ta > k \fa - 2n \ ; \ B ) \nn\\
&\le&
\Q(\ta > k \fa - 2n ) \lc \P( B ) + \phi(n)  \rc \ .
\eeqan
Apply this procedure to
$\ll \ta> (k-1)\fa\rl$ and $B=\ll \ta\circ T^{(k-1)\fa} > \fa-2n\rl$ to bound
$\Q( \ta > k \fa - 2n ) $ by
\[
    \Q( \ta > (k-1)\fa-2n )
\lc \P(\ta > \fa-2n ) + \phi(n) \rc \ .
\]
Iterate this procedure to bound $\Q( \ta > k \fa - 2n) $ by
\[\Q(\ta>\fa-2n)
\lc\P(\ta>\fa-2n)+\phi(n)\rc^{k-1} \ .
\]
This ends the proof of the Lemma.
\epf

The next proposition
establishes a relationship between {\em hitting} and {\em return}
times with an error \emph{uniform} in $t$. In particular,
(b)  says that
they are close (up to $2\eps(A)$) if and only if $\za$ is close to 1.

\bpro \label{:lambda}
Let  \pros be a \fmix processes. Let $A \in {\mathcal C}_n$ and $k$ a positive integer.
Then the following holds:
\begin{itemize}
\item[(a)] For all $0\le r \le \fa$,
\beqa
&&\lv \Q(\ta> k \fa +r) - \Q(\ta>k \fa) \P(\ta>r) \rv \\
&\le& 2 \eps(A) \ \Q(\ta> k\fa-2n)  \ .
\eeqa

\item[(b)] For all $i \ge \oa \in \Nset$,
\beqn \label{itemb}
\lv \Q(\ta>i) - \za \P(\ta>i) \rv \le 2 \eps(A) \ .
\eeqn
\end{itemize}
\epro

\bpf
To simplify notation,
for $t\in Z$ we write $\ta^{[t]}$ to mean $\ta \circ T^{t}$.
Assume $r \ge 2n$
We introduce a gap of length $2n$
after coordinate $t$ to construct the following triangule
inequality
\beqan
&&
| \Q(\ta>k\fa+r) - \Q(\ta>k\fa) \P(\ta>r)| \nn\\
&\le&\label{mix1}
|\Q(\ta>k\fa+r) - \Q(\ta> k\fa ; \ta^{[k\fa+2n]}>r-2n)|\\
&+&
|\Q(\ta>k\fa;\ta^{[k\fa+2n]}>r-2n)- \Q(\ta>k\fa) \P(\ta>r-2n)|\nn\\
   \label{mix2}\\
&+& \label{mix3}
\Q(\ta>k\fa) | \P(\ta>r-2n) -  \P(\ta>r) | \ .
\eeqan

(\ref{mix1}) is bounded by a direct computation by
$ \Q(\ta>k\fa ; \ta^{[k\fa]} \le 2n)$.
This last quantity is bounded using (\ref{1itera}) by
\[
\Q(\ta> k\fa -2n ) \lc 2n\pa+\phi(n) \rc \ .
\]

Term (\ref{mix2}) is bounded using the \fmix property by
$$ \Q(\ta>k\fa) \phi(n)  \ .$$
The modulus in (\ref{mix3}) is bounded using stationarity by
$$\P(\ta\le 2n) \le 2n \pa \ .$$
If $r<2n$, just change $r-2n$ by zero and the same proof holds.
This ends the proof of (a).

\vs              
The proof of (b) is very similar to that previous one.
We do it briefly. Write the following triangle inequality
\beqa
&&
| \Q(\ta>i) - \za \P(\ta>i)| \nn\\
&\le&
|\Q(\ta>i) - \Q(\ta> \oa ; \ta^{[\oa+2n]}>i-\oa-2n)|\\
&+&
|\Q(\ta>\oa;\ta^{[\oa+2n]}>i-\oa-2n)- \za \P(\ta>i-\oa-2n)|\nn\\
&+& 
\za | \P(\ta>i-\oa-2n) -  \P(\ta>i) | \ .
\eeqa
The moduli on the right hand side of the above inequality are bounded as follows.
The first one by 
$ \Q(\ta>\oa ; \ta^{[\oa]} \le \oa+2n-1) $
which is bounded by
$
\Q( \ta\in \ra \cup \{n,\dots,\oa+2n-1\}  ) \ .
$
The cardinal of $\ra  \cup \{n,\dots,\ta+2n-1\}$ is less or equal than $n+\oa+\ra$.
Therefore,
the last expression is bounded following the proof of Lemma \ref{lem:b1} by
$
(2n+\oa)\P(A^{(w)})+\phi(\na-w) \ .
$

The second one is bounded using the \fmix property by
$ \za \phi(n)  \ .$

The third one is bounded using stationarity by
$$\P(\ta\le \oa+2n) \le (\oa+2n) \pa .$$

This ends the proof of  (b).
\epf

The following proposition is the key of the proof of Theorem \ref{teo:ret}.

\bpro \label{indep}
Let  \pros be a  \fmix  process. 
Let $A\in\calc^n, n\in\Nset$ and let $k$ be any integer $k\ge1$.
Then the following inequality holds:
\beqa
&&
\lv \Q(\ta > k \fa) -
    \Q(\ta >   \fa) \P(\ta > \fa)^{k-1} \rv \\
&\le&
     2 \eps(A)(k-1) \Q(\ta > \fa-2n )[ \P(\ta > \fa-2n ) + \phi(n) ]^{k-2} \ .
\eeqa
\epro

\bpf 
For $k=1$ there is nothing to prove.
Take $k\ge 2$. The left hand side of the above inequality  is bounded by
\[
\sum_{j=2}^{k}
|\Q(\ta>j \fa)-\Q(\ta>(j-1)\fa)
\P(\ta>\fa)| \P(\ta>\fa)^{k-j} .
\]
The modulus in the above sum is bounded by
\[ 2 \eps(A) \Q(\ta > (j-1)\fa-2n ) \ , \]
due to Proposition \ref{:lambda} (a).
The right-most factor is bounded using Lemma \ref{itera} by
$[ \P( \ta > \fa-2n ) + \phi(n) ]^{j-2} .$
The conclusion follows.
\epf

\subsection{Proofs of Theorem  \ref{teo:ret} and corollaries }

\bpf{\bf of Theorem \ref{teo:ret} }
We divide the proof according to the different values of $t$:
(i) $t< \oa$, (ii) $\oa\le t\le \fa$ and (ii) $t > \fa$.
\vs
Consider first $t<\oa$. (\ref{imagen}) says that the left hand side of (\ref{eqteo:ret}) is zero.
\vs
Consider now $\oa \le t\le \fa$.
First write
\beqn \label{return}
\Q(\ta  > t)
= \frac{\Q(\ta  > t)}{\P(\ta>t)}\P(\ta  > t)
= p_{t+1} \P(\ta  > t) \ ,
\eeqn
and
\beqan
\P(\ta  > t)
&=&\prod_{i=\oa+1}^{t} \P(\ta>i | \ta  > i-1)  \label{hitting}\\
&=&\prod_{i=\oa+1}^{t} (1-\P(T^{-i}A | \ta  > i-1 ) ) \nn\\
&=&\prod_{i=\oa+1}^{t} ( 1- p_i \pa ) \ ,\nn
\eeqan
where
\[
p_i
\stackrel{def}{=}
\frac{\Q(\ta  > i-1 ) }{ \P(\ta>i-1) } \ .
\]
Further
\beqn \label{ineq}
\lv 1-p_i\pa -  e^{-\za\pa} \rv
\le
\lv p_i- \za \rv \pa + \lv 1-\za\pa - e^{-\za\pa} \rv \ .
\eeqn
Firstly, by Proposition \ref{:lambda} (b) and the
fact that $\P(\ta>i)\ge 1/2$ since $i\le \fa=1/(2\pa)$ we have
\beqn \label{coef}
\lv p_i - \za \rv \le {2 \eps(A) \over \P(\ta > i) }
\le 4 \eps(A) \ .
\eeqn
Secondly, note that $|1- x - e^{-x} | \le x^2/2$ for all $0 \le x \le 1$.
Apply it with $x=\za\pa$ to bound
the most right term of (\ref{ineq}) by $(\za\pa)^2/2$.
Collecting the last two bounds we get
\[
| 1-p_i\pa -  e^{-\za\pa} |
\le \frac{9}{2}  \eps(A) \pa \ , \qquad\forall i=\oa+1,\dots,\fa \ .
\]
Furthermore, since 
\beqn
|\prod a_i -\prod b_i|\le \max |a_i-b_i| (\# i)\max\{a_i;b_i\}^{\#i-1}  \qquad \forall 0 \le a_i,b_i \le 1 \ ,
\label{prodbasic}
\eeqn
we conclude from (\ref{hitting}) and (\ref{return}) that 
\beqn
|\P(\ta  > t)- e^{-\za\pa (t-\oa)}|
\le  \frac{9}{2}  \eps(A) \pa t  \ , \label{cortos1}
\eeqn
and
\beqn
|\Q(\ta  > t)- \za e^{-\za\pa (t-\oa)}|
\le  \frac{9}{2}  \eps(A) \pa t  \ , \label{cortos2}
\eeqn
for all $\ta \le t \le \fa$.
This concludes this case.

\vs
Consider now $t>\fa$.
Write it as $t=k\fa+r$ with $k$ a positive integer and $0\le r < \fa$.
We do the following triangle inequality
\beqan
&&
|\Q(\ta>t)- \za e^{-\za\pa(t-\oa)}|  \nn\\
&\le&
|\Q(\ta> k \fa +r) - \Q(\ta> k \fa)\P(\ta>r)| \label{dos}\\
&+&
|\Q(\ta>k\fa) - \Q(\ta>\fa)\P(\ta>\fa)^{k-1}  | \P(\ta>r)\label{factor1}\\
&+&
|\Q(\ta>\fa)\P(\ta>\fa)^{k-1} -\za e^{-\za k/2}| \P(\ta>r)  \label{factoriz}  \\
&+&
\za e^{-\za k/2} \ | \P(\ta>r)   - e^{-\za\pa(r-\oa)} | \label{resto}
\eeqan

By Proposition \ref{:lambda} (a), the modulus in (\ref{dos}) is bounded by
\[
2\eps(A)\Q(\ta> k \fa)  \ ,
\]
and by Lemma 4.3
\[
2\eps(A)\Q(\ta> k \fa-2n) \le 2\eps(A)(\P(\ta>\fa-2n)+\phi(n))^{k-1} \ .
\]

The modulus in (\ref{factor1}) is bounded using Proposition \ref{indep} by
\[
 2\eps(A) (k-1) (\P(\ta> \fa-2n)+\phi(n) )^{k-2} \ .
\]
Thus, the sum of  (\ref{dos}) and (\ref{factor1}) is bounded by
\beqn \label{geral}
2\eps(A)(\P(\ta> \fa-2n)+\phi(n) )^{k-2}
[ k+\phi(n) ] \ .
\eeqn
On one hand $k+\phi(n)\le k+1 \le 2k$.
On the other hand,
applying (\ref{cortos1}) with $t=\fa-2n$ we get
\beqa
|\P(\ta> \fa-2n)-e^{-\za/2+(2n+\oa)\pa }|\le \frac{9}{4}\eps(A) \ .
\eeqa
Furthermore, by the Mean Value Theorem (MVT)  
we get
\beqa
&&
|e^{-\za/2+(2n+\oa)\pa }- e^{-\za/2 }|\le 
(2n+\oa)\pa e^{(2n+\oa)\pa}  \ .
\eeqa
We conclude that for large enough $n$ 
\beq
|\P(\ta>\fa-2n)+\phi(n)-e^{-\za/2}|\le 4 \eps(A)
\eeq
And therefore (\ref{geral}) is bounded by
\beqn \label{modelo}
4\eps(A) k  (e^{-\za/2}+4\eps(A))^{k-2} \ .
\eeqn
A direct computation using Taylor's expansion gives 
\[
e^{-\za/2} \le e^{-\za/2}+4\eps(A) \le e^{-(\za/2-8\eps(A))} \ .
\]
Since $t=(k/2\pa)+r$ we get
\beq
e^{-(\za/2)(k-2)}
=e^{-\za\pa t+ \za(\pa r+1)  } \ ,
\eeq
which is bounded by
\[
e^{-\za\pa t +3/2} \ .
\]
Similarly
\beq
e^{-(\za/2-8\eps(A))(k-2)}
=e^{-(\za-16\eps(A))\pa t+ (\za-16\eps(A))(\pa r+1) } \ ,
\eeq
which for  large enough $n$ is bounded by
\[
e^{-(\za-16\eps(A))\pa t +3/2} \ .
\]
Thus  (\ref{modelo}) is bounded by
\[
36\eps(A)\pa t   e^{-(\za-16\eps(A)) \pa t }
\]

To bound (\ref{factoriz}) we proceed as follows.
From (\ref{cortos1}) and (\ref{cortos2}) with $t=\fa$ 
we get that
\beqa
&&
|\P(\ta  > \fa)- e^{-\za/2}|\\
&\le&  |\P(\ta>\fa)- e^{-\za\pa(\fa-\oa)} |+ e^{-\za/2} |e^{\za\pa\oa}-1 | \\
&\le& \frac{9}{4} \eps(A)+n\pa \\
&\le&  3 \eps(A ) \ ,
\eeqa
and similarly 
\beq
|\Q(\ta  > \fa)- \za e^{-\za/2}| \le  3  \eps(A)  \ .
\eeq
Applying the last two inequalities together with (\ref{prodbasic}), 
we get that the modulus in (\ref{factoriz}) is bounded by
\beqa
&&
3\eps(A) \ k \max\{\Q(\ta>\fa);\P(\ta>\fa); e^{-\za/2}\}^{k-1} \\
&\le& 3\eps(A) \ k 
\lp e^{-\za/2}+ 3\eps(A)  \rp^{k-1}\ .
\eeqa

An argument similar to that used to bound (\ref{modelo})
let us conclude that the last expression is bounded by
\[
10\eps(A)\pa t   e^{-(\za-12\eps(A)) \pa t } \ .
\]

The modulus in (\ref{resto}) is bounded using again (\ref{cortos1}) when $r \ge \oa$ by 
$(9/2) \eps(A)$.
If $ r<\oa$ then it can be rewritten as
\[
e^{-\za\pa (r-\oa)} -1 + \P(\ta\le r) \ ,
\]
which is bounded by
$
2n\pa .
$
We conclude that  (\ref{resto}) is bounded by
\[ (9/2) \eps(A)  e^{-\za k/2}
= (9/2) \eps(A) e^{-\za \pa (t-r)}
\le 8 \eps(A) \pa t e^{-\za \pa t} \ .
\]
This ends the proof of the theorem.
\epf

\bpf{\bf of Corollary \ref{:kac} }
Let $Y$ be the r.v. with distribution given by
\[
P(Y>t)=
\ll
\begin{array}{ll}
\ 1  &\quad \pa <t\le\pa\oa \\
\za e^{- \za(t-\pa\oa)} &\quad t<\pa\oa
\end{array} \right. \ .
\]
Then we can rewrite (\ref{eqteo:ret}) as
\begin{eqnarray}
\lv \Q(\pa\ta >t)
- \P(\pa\ta>t) \rv 
\le  C_1  \eps(A)  f(A,t/\pa) \ .  \label{int} 
\end{eqnarray}

Integrating (\ref{int}) we get
\beqan
&& \lv \bbe_A \lp (\pa\ta)^\beta \rp  - \bbe \lp Y^\beta \rp  \rv \nn \\
&=& \lv \int_{\pa}^{\infty} \beta t^{\beta-1}  \lp \P \lp \pa\ta>t \rp - \P \lp Y>t \rp \rp \rv \nn \\
&\le&   \int_{\pa}^{\infty} \beta t^{\beta-1}  \lv \P \lp \pa\ta>t \rp - \P \lp Y>t \rp     \rv \nn \\
&\le&   C_1  \eps(A) \int_{\pa}^{\infty} \beta t^{\beta-1}   f(A,t/\pa) dt \ . \nn 
\eeqan
Now we compute $\bbe \lp Y^\beta \rp =
\int_{\pa}^{\infty} \beta t^{\beta-1} \P(Y>t)  $.
We do it in each interval $[\pa,\pa\oa]$ and $[\pa\oa,\infty)$.


The first one is $(\pa\oa)^{\beta}-\pa^{\beta} $.
The second one can be re-written as 
\beqn \label{sec}
\za \  e^{ \za\pa\ta} \lp  \int_{0}^{\infty} - \int_{0}^{\pa\ta}\rp \beta t^{\beta-1}e^{- \za t} dt \ .
\eeqn
Consider the exponent of the second factor in (\ref{sec}).
By definition we have $ \za\pa\ta\le\pa n $.
Moreover, $\pa$ decays exponentially fast on $n$.
Then for the second factor we have  $ |e^{ \za\pa\ta}-1|\le C \pa n $.
Further, the first integral is ${\Gamma(\beta+1)/ \za^{\beta} }$.
The second one is bounded by $(\pa\ta)^\beta$. 
We recall that the first factor in (\ref{sec}) is $\za$.
We conclude that 
\[ 
\lv\bbe \lp Y^\beta \rp - \frac{\Gamma(\beta+1) }{ \za^{\beta-1}} \rv 
\le  Cn\pa + 2(n\pa)^\beta)
\le  C (n\pa)^{(\beta\land 1)} \ .
 \]
Similar computations give
\beqa
\int_{\pa}^{\infty} \beta t^{\beta-1} f\lp A,{t/\pa} \rp dt
&\le& 
{\beta\over\beta+1} {\Gamma(\beta+2)\over (\za-\eps(A))^{\beta+1} } \\
&\le& 
{\beta  e^{2\eps(A)(\beta+1)/\za}\over \za^2}   {\Gamma(\beta+1)\over \za^{\beta-1} } 
\ .
\eeqa
In the last inequality we used $x \le 2(1-e^{-x})$ for small enough $x>0$.
This ends the proof of the corollary.
\epf

\vs
\noindent{\bf Proof of Corollary \ref{cor:equiv}. }
$(a) \Goback (d)$. It follows directly from Theorem \ref{teo:ret}.

$(b)  \Go (a), (c)$. It follows by Theorem \ref{teo:ret} and Theorem 1 in \cite{aba4}

$(a) \Go (b)$ and $(c) \Go (b)$.
They follow by Theorem \ref{teo:ret}, Theorem 1 in \cite{aba4} and  (\ref{coef}).
The corollary is proved.
\epf

\section{SOJOURN TIME}

In this section we consider the number of consecutive visits to a fixed string $A$
and prove that the distribution law of this number can be well approximated
by a geometric law.

\bdefi
Let $A \in \calc^n$.
We define the {\em sojourn time} on the set $A$ as the r.v.
$ S_A : \Omega  \go  \Nset \cup \{\infty\}$
\[
S_A(x) = \sup \ll k \in \Nset \ | \ 
x \in A  \cap   T^{-j \oa}A \  ; \  \forall j=1, \dots , k \rl \ ,
\]
and $S_A(x) = 0$ if the supremum is taken over the empty set.
\edefi

Before to state our main result we have to introduce the following 
definition about certain continuity property of the probability $\P$
conditioned to $i$  consecutive occurrences of the string $A$.

\bdefi For each fixed $A \in \calc^n$, we define the  sequence
of probabilities $\lp \rho_i(A) \rp_{i\in\Nset}$ as follows:
\[
\rho_i(A) \stackrel{def}{=}
\P\lp A \ \big| \ \bigcap_{j=1}^{i} T^{j \oa}A \rp \ .
\]
If the limit $\lim_{n\go\infty} \rho_i(A)$ exists then we denote it by $\rho(A)$.
\edefi

\brem
By stationarity $\rho_1(A)=1-\za$.
\erem

In the following 2 examples, the sequence
$(\rho_i(A))_{i\in\Nset}$ not just converges but even is constant.

\beje \label{iid}
For a i.i.d. Bernoulli process with parameter
$0<\theta=\P(X_i=1)=1-\P(X_i=0)$, and for the $n$-string
$A=\{X_0^{n-1}=1\}$, we have that $\rho_i(A)=1-\za=\theta$ for all $i\in\Nset$.
\eeje

\beje \label{markov}
Let \pros be a  irreducible and aperiodic finite state Markov chain.
For $A=\{X_0^{n-1}=a_{0}^{{n-1}} \}\in \calc^n$,
the  sequence $ \lp \rho_i(A) \rp_{i\in\Nset}$  is constant.
More  precisely,  by the  Markovian property and for all $i \in \Nset$
\beqa
\rho_i(A)
&=&\P\lp X_{n-\oa}^{n-1}= a_{n-\oa}^{n-1}   |  X_{n-\oa-1}= a_{{n-\oa-1}} \rp \\
&=&\prod_{j=n-\oa}^{n-1} \P\lp X_{j}= a_{j} | X_{j-1}= a_{j-1} \rp
\ .
\eeqa
\eeje

The next is an example of a process with infinity memory and converging  
$(\rho_i(A))_{n\in\Nset}$.

\beje \label{esparça}
The following is a family of processes of the renewal type.
Define $(X_n)_{n\in\Nset}$ as the order one Markov chain over $\Nset$
with transitions probabilities given by
\[
Q(n,n+1)=q_n   \qquad Q(n,0)=1-q_n \qquad \forall n\ge 0
\]
Define the process
\[ Y_n = \ll
\begin{array}{ll}
0 & X_n=0 \\
1&  X_n\not=0
\end{array} \right.
\]
The process $(X_n)_{n\in\Nset}$ is positive recurrent (and  then $(Y_n)_{n\in\Nset}$)
if and only if 
$\sum_{k=0}^{\infty} \prod_{j=0}^{k} q_j < \infty $.
Direct computations show that 
$$P(Y_{0}^{n-1}=1)= \sum_{k=0}^{n} \prod_{j=0}^{k} q_j \qquad \forall n\in\Nset \ .$$
Now chose $q_j$ such that $P(Y_{0}^{n-1}=1)= e^{-n+\del(n)}$
with $\del(n)$ any converging sequence (to any real number)
and such that $|\delta(i+1)-\delta(i)|< 1$ for all $i\in\Nset$.
Take $A=\{Y_{0}^{n-1}=1\}$. Thus $\oa=1$.
Then 
\[
\rho_i(A)=e^{-1+\del(n+i+1)-\del(n+i)} \qquad  {\rm and } \qquad
\lim_{i\go\infty} \rho_i(A)=e^{-1} \in (0,1) \ . 
\]
\eeje

\vs
In the following theorem we assume that $(\rho_i(A))_{i \in \Nset}$ converges
with velocity $d_i(A)$. Namely, there is a real number $\lima\in[0,1)$
such that
\beqn
\lv \rho_i(A) - \lima \rv \le d_i(A) \qquad
{\rm for\ all\ } i \in \Nset,  \label{eq:lim}
\eeqn
where $d_i$ is a sequence converging to zero.
For simplicity we put $\sd(A)=\sup\{d_i(A) \ | \ i\in \Nset\}$.

\bteo \label{teo:geo}
Let \pros be a stationary process. Let $A \in \calc^n$.
Assume that (\ref{eq:lim}) holds.
Then, there is $c(A)\in [0,1)$, such that the following
inequalities hold for all $k \in \Nset$:
\beq
\lv \Q \lp S_A = k \rp - (1-\lima)\ \lima^k \rv
\le   c(A)^{k} \sum_{i=1}^{k+1} d_i(A)
\le   c(A)^{k} (k+1) \sd(A)  \ .
\eeq
\eteo

We deduce immediately that the $\beta$-moments
of $S_{A}$ can be approximated by $\bbe(Y^\beta)$  where 
$Y$ is a geometric random variable with parameter $\rho(A)$.

\bcor \label{:mome} Let $Y$ be a r.v. with geometric distribution with
parameter $\lima$.
Let $\beta >0$.
Then
\beq
\lv \bbe_A  \lp S_A^\beta \rp  -  \bbe(Y^\beta)\rv
\le 2 \sd(A) \sum_{k=1}^{\infty} k^{\beta+1} c(A)^{k} \ .
\eeq
\ecor

\brem The sum $\sum_{k=1}^{\infty} k^{\beta+1} c(A)^{k}$
can be approximated using the Gamma function by $\Gamma(\beta+2)/ (-\ln c(A))^{\beta+2}$.
When the supremum of the distances $|\rho_i(A)-\rho(A)|$ is small, 
the approximations given by Theorem \ref{teo:geo} and Corollary \ref{:mome} 
are good. The smaller is $c(A)$, the better they are.
We compute these quantities for the examples of this section.
\erem

\noindent{\bf Example \ref{iid}} (continuation) 
It follows straight-forward from definitions that 
$\rho_i(A)=\rho(A)=\P(A^{(\oa)})$ for all $i$ and for any $A\in \calc^n, n\in \Nset$.
Thus $c(A)=\P(A^{(\oa)})$ and $\sd(A)=0$.
\vs
\noindent{\bf Example \ref{markov}} (continuation) 
We already compute that
$\rho_i(A)=\rho(A)$ for all $i$ and for any $A\in \calc^n, n\in \Nset$.
Thus $c(A)=\rho(A)$ and $\sd(A)=0$.
\vs
\noindent{\bf Example \ref{esparça} }(continuation) 
For the same $n$-string there considered, we have
\[
d_i(A) = e^{-1}|e^{\del(n+i+1)-\del(n+i)}-1|
\le | \del(n+i+1)-\del(n+i) | \ ,
\]
and
\[
\sd(A) \le  \sup\{ |\del(n+i+1)-\del(n+i)| , i\in\Nset\} \ .
\]
So, for large enough $n$, $\sd(A)$ is small.
Finally, 
\[
c(A) = \max\{ e^{-1}; \sup_{n\in\Nset} e^{-1+\del(n+i+1)-\del(n+i)}  \} \in (0,1) \ .
\]
\vs
In the proof of Theorem \ref{teo:geo} we will use the following lemma. 

\blem \label{lem:prod}
Let $(l_i)_{i\in\Nset}$ be a sequence of real numbers such that
$0\le l_i<1$, for all $i \in \Nset$. Let $0\le l<1$ be such that
$|l_i - l| \le d_i$ for all $i \in \Nset$ with $d_i \go 0$.
Then, there is a constant $c \in [0,1)$, such that the
following inequalities hold for all $k \in \Nset$:
\[
\lv \prod_{i=1}^{k} l_i - l^k \rv
\le c^{k-1 }\sum_{i=1}^{k} d_{i}
\le k \ c^{k-1} \sd \ .
\]
where $\sd=\sup\{d_i, i\in\Nset\}$.
\elem

\bpf
\beqan
\lv \prod_{i=1}^{k} l_i - l^k \rv
&=&
\lv \prod_{i=1}^{k} l_i - \prod_{i=1}^{k-1} l_i l  + \prod_{i=1}^{k-1} l_i l
                        - \prod_{i=1}^{k-2} l_i l^2+ \prod_{i=1}^{k-2} l_i l^2
- \dots - l^k \rv \nn \\
&\le& \sum_{i=1}^{k} \lp \prod_{j=1}^{k-i} l_j \rp  |l_{k-i+1} -l| \ l^{i-1}
 \le  c^{k-1 }\sum_{i=1}^{k} d_{i}  \nn \\
&\le& k \  c^{k-1 } \sd \nn \ ,
\eeqan
where $c = \max\ll l ; \sup_{i\in\Nset} l_i  \rl$.
\epf

\bpf{\bf of Theorem \ref{teo:geo} }
For  $k=0$, we just note that $\Q \lp S_A = 0 \rp = 1-\rho_1(A) $ and
$ |1- \rho_1(A) - (1-\lima) | \le d_1(A)$.
Suppose  $k \ge 1$. Therefore
\beqa
&& \Q\lp S_A = k \rp \\
&=&\Q\lp   \bigcap_{j=0}^{k}T^{-j \oa}A  \ \ ; \ \ T^{-(k+1)\oa}A^c    \rp\\
&=&\P\lp T^{-(k+1)\oa}A^c  |  \bigcap_{j=0}^{k} T^{-j \oa}A \rp
   \prod_{i=1}^{k}\P\lp T^{-i \oa}A  | \bigcap_{j=0}^{i-1} T^{-j \oa}A \rp\\
&=&(1-\rho_{k+1}(A))\prod_{i=1}^{k} \rho_i(A) \ .
\eeqa
Third equality follows by stationarity.
Lemma \ref{lem:prod} ends the proof of the theorem.
\epf

\bpf{\bf of Corollary \ref{:mome}}
We use the inequality
\[
\lv \bbe \lp X^\beta \rp  - \bbe \lp Y^\beta \rp  \rv
\le
\sum_{k\ge 0} k^\beta
\lv \P \lp X=k \rp - \P \lp Y=k \rp \rv
\ ,
\]
which holds for any pair of positive r.v. $X, Y$.
We apply the above inequality with $X=S_A$ and $Y$ geometrically distributed
with parameter $\rho(A)$.

The exponential decay of the error term in
Theorem \ref{teo:geo} ends the proof of the corollary.
\epf

{\bf Acknowledgments}
The authors are beneficiaries of a Capes-Cofecub grant.
We thank P. Ferrari and  A. Galves for useful discussions.
We kindly thank also two anonymous referees for their useful comments and suggestions
to improve a previous version of this article.


\begin{thebibliography}{10}


\bibitem{aba1}
Abadi, M. (2001).
 Exponential approximation for hitting times in
          mixing processes.
{\em Math. \ Phys. \ Elec. J. \bf 7}, 2.

\bibitem{aba4}
 Abadi, M. (2004).
 Sharp error terms and necessary conditions for exponential
hitting times in mixing processes.
 {\em Ann. \ Probab. \bf 32}, 1A, 243-264.

\bibitem{AG}
 Abadi, M., and Galves, A. (2001).
 Inequalities for the occurrence times of rare events in mixing
          processes. The state of the art .
{\sl Markov \ Proc. \ Relat. \ Fields. \bf 7}, 1, (2001) 97-112.

\bibitem{AHV}
 Abadi, M., and Vaienti, S. (2006).
 Statistics properties of repetition times.
{\sl Preprint.} 

\bibitem{BT}
Blom, G., and Thorburn D. (1982).
How many random digits are required until given sequences are obtained?
{\sl J. App. Prob.\bf 19} 518-531.

\bibitem{bow}
 Bowen, R. (1975).
 Equilibrium states and the ergodic theory of
          Anosov diffeomorphisms.
 {\sl Lecture Notes in Math,  \bf 470.
\sl Springer-Verlag, New York.} 

\bibitem{chaz}
 Chazottes, J.-R. (2003).
 Hitting and returning to non-rare events in mixing dynamical systems.
 {\sl Nonlinearity \bf 16}, 1017-1034.


\bibitem{CFS}
 Cornfeld, I., Fomin, S., and Sinai Y. (1982).
 Ergodic theory.
 {\sl Grundlähren der Mathematischen Wissenschaften, \bf 245.
\sl Springer-Verlag, \ New York.} 


\bibitem{dou}
 Doukhan, P. (1995).
 Mixing. Properties and examples.
 {\sl Lecture \ Notes \ in \ Statistics \bf 85,
\sl Springer-Verlag.} 

\bibitem{GS}
 Galves, A., and Schmitt, B. (1997)
 Inequalities for hitting times in mixing dynamical systems.
 {\sl Random \ Comput. \ Dyn. \bf 5}, 337-348.

\bibitem{hir}
 Hirata, M. (1993).
 Poisson law for Axiom A diffeomorphism.
 {\sl Ergod. Th. Dyn. Sys. \bf 13}, 533-556.

\bibitem{HSV}
 Hirata, M., Saussol, B. and  Vaienti, S. (1999).
 Statistics of return times: a general framework and new applications.
 {\sl Comm. \ Math. \ Phys. \bf 206}, 33-55.

\bibitem{kac}
 Kac, M. (1947).
 On the notion of recurrence in discrete stochastic processes.
 {\sl Bull.\ Amer.\ Math.\ Soc. \bf 53}, 1002-1010.

\bibitem{LSV}
Liverani, C., Saussol, B. and  Vaienti, S. (1998).
Conformal measures and decay of correlations for covering weighted 
systems.
{\sl Ergod. Theeor. dynam. Sys. \bf 18} 1399-420.

\bibitem{OW}
Ornstein, D., and Weiss, B. (1993).
Entropy and data compression schemes.
{\sl IEEE \ Trans. \ Inform. \ Theory \bf 39}, 1,  78-83.

\bibitem{ste}
Stefanov, V. (2003).
The intersite distances between pattern occurrences in strings generated by
general discrete- and continuous-time models. An algorithmic approach
{\sl J. App. Prob. \bf 40}, 881-892. 


\bibitem{wal}
Walters, P. (1975).
Ruelle's operator theorem and $g$-measures.
{\sl Trans. Amer. Math. Soc. \bf 214}, 375-387.

\bibitem{WZ}
Wyner, A., and  Ziv, J. (1989).
Some asymptotic properties of the entropy of a stationary
ergodic data source with applications to data compression.
{\sl IEEE \ Trans. \ Inform. \ Theory \bf 35}, 6, 1250-1258.

\end{thebibliography}
\end{document}